\documentclass[12pt]{amsart}
\usepackage{amssymb,mathrsfs,longtable}
\usepackage[matrix,arrow,curve]{xy}
\usepackage{setspace}
\usepackage{multirow}
\usepackage{textcomp}
\usepackage{hyperref}

\usepackage[dvipsnames]{xcolor}
\usepackage{graphicx}

\usepackage[normalem]{ulem}
\usepackage{cancel}

\newcounter{cequation}[section]

\newtheorem{theorem}[cequation]{Theorem}
\newtheorem*{theorem*}{Theorem}
\newtheorem{lemma}[cequation]{Lemma}
\newtheorem{corollary}[cequation]{Corollary}

\newtheorem{proposition}[cequation]{Proposition}

\theoremstyle{definition}

\newtheorem{definition}[cequation]{Definition}
\newtheorem*{definition*}{Definition}

\newtheorem{problem}[cequation]{Problem}
\newtheorem*{notation*}{Notation}

\theoremstyle{remark}
\newtheorem{remark}[cequation]{Remark}

\makeatletter\@addtoreset{equation}{section}
\makeatletter\@addtoreset{section}{part}

\makeatother

\def \O {\mathcal{O}}
\def \CC {\mathbb{C}}
\def \P {\mathbb{P}}
\def \PP {\mathbb{P}}

\def \QQ {\mathbb{Q}}

\def \RR {\mathbb{R}}

\def \Pic {\mathrm{Pic}}

\def \RRR {\mathcal{R}}
\def \hh {\mathrm{h}}

\def \ge {\geqslant}
\def \le {\leqslant}

\title{On Hodge level of weighted complete intersections of general type}

\author{Victor Przyjalkowski} 

\address{\emph{Victor Przyjalkowski}
\newline
\textnormal{Steklov Mathematical Institute of RAS, 8 Gubkina street, Moscow 119991, Russia.
}
\newline
\textnormal{\texttt{victorprz@mi-ras.ru, victorprz@gmail.com}}}
\thanks{
This work was performed at the Steklov International Mathematical Center and supported by the Ministry of Science and Higher Education of the Russian Federation (agreement no. 075-15-2022-265).
The author was supported by the Foundation for the
Advancement of Theoretical Physics and Mathematics ``BASIS''}

\begin{document}

\begin{abstract}
We show that smooth varieties of general type which are well formed weighted complete intersections of Cartier divisors have maximal Hodge level,
that is, their the rightmost middle Hodge numbers do not vanish.
We show that this does not hold in the quasi-smooth case.
\end{abstract}

\maketitle

\section{Introduction}

One of the most basic biregular invariants of smooth varieties are Hodge numbers.
In particular, if a lot of them vanish, then the variety tends to be ``simple'' from
homological point of view. In  \cite[\S1]{Rapoport}, \cite[\S2a]{Carlson}, and~\cite{PSh20a}
for a smooth variety $X$ it was introduced the notion of \emph{Hodge level}
$\hh(X)=\max\{q-p\mid h^{p,q}(X)\neq 0\}$. In other words, this number measures
how big the number $a$ such that $HH_a(D^b(Coh X))\neq 0$ can be.

Let $n=\dim (X)$. By definition, $\hh(X)$ has maximal possible value $n$
if and only if~$h^{0,n}(X)\neq 0$.
If $X$ is Fano, that is, if its anticanonical class is ample,
then one has~$h^{0,n}(X)=h^0(\Omega_X^n)=0$. If $X$ is Calabi--Yau,
then, by definition, $H^{0,n}(X)=\CC$, so it has maximal possible Hodge level.
Obviously, the Hodge level of curves of positive genus is equal to $1$, so it is maximally possible again.
However already in dimension $2$ there exist examples of \emph{fake projective planes} --- smooth surfaces of general type (that is, whose canonical class is nef and big), such that their
Hodge diamond is the same as for $\PP^2$, see, for example,~\cite{Mu79},~\cite{CS10}.

To study the Hodge level of varieties one needs to be able to compute Hodge numbers.
However this, as well as even to determine if a variety is
Fano, Calabi--Yau or of general type (or none of them), can be not easy. The most straightforward way to construct varieties
is to describe them as complete intersections in something well studied, like
in toric varieties or homogenous spaces. There is an approach, going back to Griffiths, to compute
Hodge numbers of complete intersections in toric varieties (see~\cite{Di95},~\cite{Gr69},~\cite{Na97},~\cite{Ma99});
its generalization to some of complete intersections in Grassmannians can be found
in~\cite{FM18}. Probably the most complicated and interesting (from different points of view)
class of varieties are ones with Picard rank one.
Let $X$ be a quasi-smooth complete intersection in a toric variety $T$.
By Lefschetz-type theorem
for complete intersections in toric varieties (see~\cite[Proposition 1.4]{Ma99}) this means that
$T$ is either a weighted projective space or its quotient by finite
group, see, for instance,~\cite{Ka09}. If $X$ in addition is a smooth Fano variety, then $T$ in fact is a weighted projective space by~\cite[Theorem 2]{PSh20c},
so $X$ is a \emph{weighted complete intersection}.
This is not true if $X$ is of general type, see~\cite{PSh20c}.
However the weighted complete intersection case is still the most approachable.

Hodge levels of Fano weighted complete intersections are studied in~\cite{PSh20a}.
More precise, in loc. cit. it is studied how \emph{small} the Hodge level
of Fano weighted complete intersections can be, and the varieties of smallest or close to be smallest Hodge levels are
classified. In the present paper, on the contrary, we study \emph{how big Hodge levels
of weighted complete intersections of general type can be}.
We expect that in this case the phenomena of fake projective planes does not
appear, and Hodge levels for smooth weighted complete intersections of general type are
maximally possible.

The following summarizes the speculation above and shows that our expectation holds
in the case of complete intersections of Cartier divisors.

\begin{theorem}
\label{theorem:maximal Hodge level}
Let $X\subset\PP(a_0,\ldots,a_N)$ be a smooth well formed weighted complete intersection
of multidegree $(d_1,\ldots,d_k)$. Put $n=N-k=\dim (X)>0$ and $i_X=\sum d_u-\sum a_l$.
Then
\begin{itemize}
      \item[(i)] Let $X$ be Fano variety, that is $i_X<0$. Then $\hh(X)<n$.
      \item[(ii)] Let $X$ be Calabi--Yau variety, that is $i_X=0$. Then $\hh(X)=n$. Moreover, $h^{0,n}(X)=1$.
      \item[(iii)] Let $X$ be 
      a complete intersection of Cartier divisors, and let it be of general type, that is 
      $i_X>0$. Then $\hh(X)=n$.
    \end{itemize}
\end{theorem}

\begin{corollary}
In the terminology of~\cite{PSh20a} this, in particular, means that smooth well formed weighted Calabi--Yau complete intersection of dimension $n$
is of $n$-Calabi--Yau type, while if a smooth $n$-dimensional well formed weighted complete intersection of Cartier divisors is of general type,
then it is never $\QQ$-homologically minimal, Hodge--Tate, of curve type (if $n>1$), or of $m$-Calabi--Yau type for $m<n$.
\end{corollary}

In Corollary~\ref{corollary:codim 2 wci} we show that the assertion (iii) of Theorem~\ref{theorem:maximal Hodge level}
holds for codimension $2$ smooth weighted complete intersections, not necessary of Cartier divisors.
We expect that this holds in the general case, for all smooth weighted complete intersections of general type (Problem~\ref{problem:general type maxmal}).

Assertions (i) and (ii) of Theorem~\ref{theorem:maximal Hodge level} hold for quasi-smooth well formed weighted complete intersections
(see Corollary~\ref{corollary:h0n}). However assertion (iii) can fail even in the hypersurface case.
The example is given in Proposition~\ref{proposition:quasi-smooth hypersurface}.

Note that $\hh(X)=n$ means that $h^{0,n}(X)>0$, that is $\dim |\mathcal{O}_X(i_X)|\ge 0$.
In fact~\cite[Conjecture 4.8]{PST17} claims that $\dim |\mathcal{O}_X(m)|\ge 0$ for all $m\ge i_X$.

\medskip

In Section~\ref{section:preliminaries} we give some definitions and results related to the subject.
In Section~\ref{section:Cartier} we prove Theorem~\ref{theorem:maximal Hodge level} and Corollary~\ref{corollary:codim 2 wci}.
In Section~\ref{section:quasi-smooth} we provide the counterexample to the statement of Theorem~\ref{theorem:maximal Hodge level}
in the quasi-smooth case (Proposition~\ref{proposition:quasi-smooth hypersurface}).

\medskip

The author is grateful to M.\,Korolev for the proofs of Lemma~\ref{lemma:2n_primes} and Proposition~\ref{proposition:prime sequence},
for C.\,Shramov for helpful comments and Proposition~\ref{proposition:codim 2}, and for the referee whose useful remarks
improved the paper.

\section{Preliminaries}
\label{section:preliminaries}
Let us introduce basic definitions and results related to weighted complete intersections.
More details see in~\cite{Do82},~\cite{IF00}, and~\cite{PSh}.

Let $a_0,\ldots,a_N$ be positive integers. Consider the graded algebra~\mbox{$\CC[x_0,\ldots,x_N]$},
where the grading is defined by assigning the weights $a_l$ to the variables~$x_l$.
Put
$$
\P=\P(a_0,\ldots,a_N)=\mathrm{Proj}\,\CC[x_0,\ldots,x_N].
$$
We use the abbreviation
\begin{equation*}
(a_0^{k_0},\ldots,a_m^{k_m})=
(\underbrace{a_0,\ldots,a_0}_{k_0},\ldots,\underbrace{a_m,\ldots,a_m}_{k_m}),
\end{equation*}
where $k_0,\ldots,k_m$ are allowed to be
any positive integers.
If some of $k_i$ is equal to $1$ we drop it for simplicity.
The weighted projective space $\P$ is said to be \emph{well formed} if the greatest common divisor of any $N$ of the weights~$a_l$ is~$1$. Every weighted projective space is isomorphic to a well formed one, see~\cite[1.3.1]{Do82}.
A subvariety $X\subset \P$ is said to be \emph{well formed}
if~$\P$ is well formed and
$$
\mathrm{codim}_X \left( X\cap\mathrm{Sing}\,\P \right)\ge 2,
$$
where the dimension of the empty set is defined to be $-1$.

We say that a subvariety $X\subset\P$ of codimension $k$ is a \emph{weighted complete
intersection of multidegree $(d_1,\ldots,d_k)$} if its weighted homogeneous ideal in $\CC[x_0,\ldots,x_N]$
is generated by a regular sequence of $k$ homogeneous elements of degrees $d_1,\ldots,d_k$.
A weighted complete intersection~$X\subset\P$
is said to be \emph{an intersection
with a linear cone} if one has $d_u=a_l$ for some~$u$ and~$l$.
%
%
%
Let $p\colon \mathbb A^{N+1}\setminus \{0\}\to \P$ be the natural projection to the weighted projective space. A subvariety $X\subset \P$
is called \emph{quasi-smooth} if the preimage~\mbox{$p^{-1}(X)$} is smooth.
Note that by~\cite[Proposition 2.9]{PSh20b} if a quasi-smooth weighted complete intersection $X\subset \PP$ of dimension at least $3$ is general in the family of weighted complete intersections of the same multidegree in $\PP$, then there exists a quasi-smooth well formed weighted
complete intersection isomorphic to $X$ which is not an intersection with a linear cone.

%

Smoothness of a weighted complete intersection implies arithmetic restrictions on weights and degrees defining it.

\begin{proposition}[{cf.~\cite[Proposition~4.1]{ChenChenChen}}]
\label{proposition:smooth-condition}
Let $X\subset\PP(a_0,\ldots,a_N)$ be a smooth well formed weighted complete intersection
of multidegree $(d_1,\ldots,d_k)$.
Then for every $r$ and every choice of~$r$ weights~\mbox{$a_{i_1},\ldots,a_{i_r}$, $i_1<\ldots<i_r$},
such that their greatest common divisor $\delta$ is greater than~$1$,
there exist $r$ degrees $d_{s_1},\ldots,d_{s_r}$, $s_1<\ldots<s_r$,
such that their greatest common divisor is divisible by~$\delta$.
\end{proposition}

In a similar way one can determine if a hypersurface in a weighted projective space is a Cartier divisor.

\begin{proposition}[{\cite[Proposition 8]{RoTe12} or
the proof of~\mbox{\cite[Theorem 3.2.4(i)]{Do82}}}]
\label{proposition:Cartier}
Let~\mbox{$\PP=\PP(a_0,\ldots,a_N)$} be a well formed weighted
projective space.
Then the Picard group $\Pic(\PP)$ is a free group generated by $\O_\PP(r)$,
where $r$ is the least common multiple of the weights~$a_l$.
In particular, a degree $d$ hypersurface is Cartier if and only if
$a_l|d$ for~\mbox{$l=0,\ldots,N$}.
\end{proposition}

By combinatorial reasons it is more convenient to deal just with the collections of weights and degrees defining
weighted complete intersections, ignoring geometric objects they provide.

\begin{definition}[cf.~\cite{PST17}]
The pair $(\overline{d},\overline{a})$ of collections of positive integers
$$\overline{d}=(d_1,\ldots,d_k),\qquad\overline{a}=(a_0,\ldots,a_N)
$$
is called
a \emph{regular pair} if the divisibility conditions of Proposition~\ref{proposition:smooth-condition}
hold for them.
The numbers $d_u$ are called \emph{degrees} and the numbers $a_l$ are called \emph{weights}.
We call the regular pair \emph{Cartier} if for all $l$ and $u$ one has $a_l|d_u$.
We call the regular pair \emph{Fano} if for
$i=\sum d_u-\sum a_l$ one has $i<0$, \emph{Calabi--Yau} if $i=0$, and \emph{of general type} if $i>0$.
\end{definition}

\begin{remark}
\label{remark: regular linear cone}
Let $(\overline{d},\overline{a})$ be a regular pair such that, up to permutation, one has
$$
\overline{a}=(d_1,\ldots,d_k,a_{k},\ldots,a_{N}).
$$
Then $\overline{a}=(d_1,\ldots,d_k,1^{N+1-k})$.
Indeed, if, say, $a_k>1$, then, by regularity, one of degrees, say $d_1$, is divisible by $a_k$;
this means that, since two weights $d_1$ and $a_k$ are divisible by $a_k$, one more degree is divisible by $a_k$, and so on.
This shows that all degrees $d_1,\ldots,d_k$ are divisible by $a_k$, so $a_k$ divides at least $k+1$ weights,
which contradicts with regularity of the pair.
\end{remark}

Weights and degrees defining smooth (or, in addition, intersection of Cartier divisors, or Fano, or Calabi--Yau, or of general type)
weighted complete intersection form a regular pair (satisfying the additional requirements).
Since the proofs of the statements related to the subject are obtained by analysis of the numerical conditions
on weights and degrees (and the statements are often formulated in terms of regular pairs),
it is usually more convenient to use the language of regular pairs
and then to derive the results for weighted complete intersections from ones for regular pairs.
We also will follow this strategy.
We refer to the statements formulated for smooth weighted complete intersections
(which follow from the statements for regular pairs)
replacing the weighted complete intersections by regular pairs,
and then we derive the final statements for weighted complete intersections from them.

It 
turns out that there is a strong restriction on the minimal weight for Fano and Calabi--Yau regular pairs.


\begin{theorem}[{cf. \cite[Corollary 3.4]{PSh20b}}] 
\label{corollary:PST}
Consider a Fano or Calabi--Yau regular pair $(\overline{d},\overline{a})$, where $\overline{d}=(d_1,\ldots,d_k)$ and $\overline{a}=(a_0,\ldots,a_N)$.
Put $i=\sum d_u-\sum a_l$. Let $a_l\neq d_u$ for all $l$ and $u$.
Let $a_0\le\ldots\le a_N$.
Then
$a_{k-i-1}=1$. Moreover, $a_{k-i}=1$ unless $\overline{d}=(6^k)$ and $\overline{a}=(1^s,2^k,3^k)$.
\end{theorem}

\begin{corollary}
\label{corollary:no 1}
Consider a regular pair $(\overline{d},\overline{a})$, where $\overline{d}=(d_1,\ldots,d_k)$, $\overline{a}=(a_0,\ldots,a_N)$,
and $N\ge k$.
Assume that $a_l>1$ for all $l$.
Then $(\overline{d},\overline{a})$ is of general type.
\end{corollary}

\begin{proof}
By Remark~\ref{remark: regular linear cone} there exists a degree $d_u$ such that $d_u\neq a_l$ for all $l$.
Deleting pairs $d_u$ and $a_l$ such that $d_u=a_l$ from the regular pair, one gets a regular pair satisfying
conditions of Theorem~\ref{corollary:PST}. It remains to apply it.
\end{proof}


Now recall how to compute Hodge numbers of weighted complete intersections.
Let $X\subset \PP=\PP(a_0,\ldots,a_N)$ be a quasi-smooth weighted complete
intersection of hypersurfaces of degrees $d_1,\ldots,d_k$.
By~\cite[\S3]{Ma99}
or~\mbox{\cite[\S11]{BC94}}, there is
a pure Hodge structure on the cohomology of $X$.
In particular, Hodge numbers~\mbox{$h^{p,q}(X)$} are well defined.
By~\cite[Theorem 10.8 and Remark 10.9]{Da78}, for the weighted projective space $\PP$
one has $h^{p,q}(\PP)=1$ if $p=q$ and $h^{p,q}(\PP)=0$ otherwise.
By~\cite[Proposition 3.2]{Ma99},
the only Hodge numbers of
~$X$
that are not inherited from
the ambient weighted projective space are $h^{p,q}(X)$ with~\mbox{$p+q=\dim(X)$}.

Put
$$
S'=\CC[x_0,\ldots,x_{N}],
$$
where the weight of $x_i$ is $a_i$, and
$$
S=\CC[x_0,\ldots,x_{N}, w_1,\ldots,w_k].
$$
Let $f_1,\ldots,f_k$ be polynomials in $S'$ of weighted degrees $d_1,\ldots,d_k$
that generate the weighted homogeneous ideal of the weighted complete intersection $X$. 
Let
$$
F=F(f_1,\ldots,f_k)=
w_1f_1+\ldots+w_kf_k\in S.
$$
Denote by $J=J(F)$ the ideal in $S$ generated by
\begin{equation*}\label{eq:generators-J}
\frac{\partial F}{\partial w_1},\ldots, \frac{\partial F}{\partial w_k},
\frac{\partial F}{\partial x_0},\ldots,\frac{\partial F}{\partial x_{N}}.
\end{equation*}
Put
$$
\RRR=\RRR(f_1,\ldots,f_k)=S/J.
$$

{The algebra $S$ is bigraded by $\deg (x_l)=(0,a_l)$ and $\deg(w_u)=(1,-d_u)$,
so that $F$ is a bihomogeneous polynomial of bidegree~$(1,0)$.
Thus, the bigrading descends to the ring~$\RRR$.}

Put $n=N-k=\dim (X)$
and
\begin{equation*}\label{eq:i-x}
i_X=\sum d_u-\sum a_l.
\end{equation*}
Let $h_{pr}^{n-q,q}(X)$ be primitive
middle Hodge numbers of~$X$, that is,
$$h_{pr}^{p,q}(X)=h^{p,q}(X)$$ for $p\neq q$ and
$$h_{pr}^{p,q}(X)=h^{p,q}(X)-1$$ otherwise.

\begin{theorem}[{see~\cite{Di95},~\cite{Gr69},~\cite[{Proposition~2.16}]{Na97},~\cite[{Theorem 3.6}]{Ma99}}]
\label{theorem:middle-Hodge-numbers}
One has
$$
h_{pr}^{q,n-q}(X)=\dim \RRR_{q,i_X}.
$$
\end{theorem}

\begin{corollary}
\label{corollary:h0n}
Denote by $J'$ the ideal in $S'$ generated by $f_1,\ldots,f_k$ and
let
$$
\RRR'=\RRR'(f_1,\ldots,f_k)=S'/J'
$$
be a graded polynomial ring.
Then $h^{0,n}=\RRR'_{i_X}$.
In particular, $h^{0,n}(X)>0$ if and only if there exists a weighted monomial of degree $i_X$ in $S'$.
\end{corollary}
\begin{proof}
One may assume that $X$ does not line in any coordinate hyperplane,
that is $x_r\notin J'$ for any $r$.
Let $M$ be a monomial of degree $i_X$ in $S'$. Then $M\notin J'$, otherwise $X$ is reducible.
Thus, $M$ descends to a non-trivial element in $\RRR'$.
On the other hand, any monomial summand of a homogenous element in $S'$ of degree $i_X$
has the same degree $i_X$.
\end{proof}

\begin{corollary}
\label{corollary:sum decomposition}
One has $h^{0,n}(X)>0$ if and only if there exist positive integers $\beta_0,\ldots,\beta_N$
such that
\begin{equation*}
\label{equation:splitting}
d_1+\ldots+d_k=\beta_0\cdot a_0+\ldots+\beta_N\cdot a_N.
\end{equation*}
\end{corollary}

\begin{proof}
If the numbers exists, then the monomial $x_0^{\beta_0-1}\cdot\ldots\cdot x_N^{\beta_N-1}$
is a weighted monomial of degree $-i_X$ in $S'$. From the other hand, if
$x_0^{\alpha_0}\cdot\ldots\cdot x_N^{\alpha_N}$ is a monomial of degree $-i_X$,
then $x_0^{\alpha_0+1}\cdot\ldots\cdot x_N^{\alpha_N+1}$ is the monomial of degree
$d_1+\ldots+d_k$, so that one can define $\beta_l=\alpha_l+1$.
\end{proof}

\begin{definition}[{\cite[Definition 1.14]{PSh20a}, cf. \cite[\S1]{Rapoport}, \cite[\S2a]{Carlson}}]
\label{definition: Hodge level}
For a smooth projective variety  $X$ put
$$
\hh(X)=\max\{q-p\mid h^{p,q}(X)\neq 0\}.
$$
The number $\hh(X)$ is called \emph{the Hodge level} of~$X$.
\end{definition}

\section{Weighted intersections of Cartier divisors of general type}
\label{section:Cartier}

In this section we prove Theorem~\ref{theorem:maximal Hodge level}. The key ingredient for this is the following.

\begin{proposition}
\label{proposition:h0n Cartier general type}
Let $(\overline{d},\overline{a})$, where $\overline{d}=(d_1,\ldots,d_k)$ and $\overline{a}=(a_0,\ldots,a_N)$, be a Cartier regular pair of general type,
and let $N-k>0$.
Then there exist positive integers $\beta_0,\ldots, \beta_N$ such that
\begin{equation}
\label{equation:splitting}
\tag{$\bigotimes$}
d_1+\ldots+d_k=\beta_0\cdot a_0+\ldots+\beta_N\cdot a_N.
\end{equation}
\end{proposition}

\begin{proof}
Let us prove this by induction on $d_1+\ldots+d_k$. The base (when $\overline {d}=4$ and $\overline{a}=(1,1,1)$) is trivial.
Denote $i=\sum d_u-\sum a_l>0$.
Let $a_l=1$ for some $l$.
Then we can choose $\beta_m=1$, $m\neq l$, and $\beta_l=i+1$.
Thus we may assume that $a_l>1$ for all $l$.
Moreover, we may assume that $a_l\neq d_u$ for all $l$ and $u$.
Indeed, otherwise we can set $\beta_l=1$ and decrease $d_1+\ldots+d_k$ removing $a_l$ and $d_u$;
clearly, after the removing we get Cartier regular pair of general type with $N-1>k-1$,
so the assertion of the proposition in this case can be given by induction.

Assume that $a_l>1$ for all $l$ and $a_l\neq d_u$ for all $l$ and $u$.
Choose a prime $p$.
Let $s$ be a maximal $p$-adic valuation for $a_0,\ldots,a_N$.
We may assume that for some integer $k\ge r\ge 1$ the numbers $a_0,\ldots a_{r-1}$ are divisible by $p^s$, while 
$a_{r},\ldots,a_{N}$ are not. Then by the divisibility assumption
$p^s|d_m$ for all $m=1,\ldots,k$.
Put
\begin{multline*}
d_u=pd_u', \ u=1,\ldots,k,\qquad a_l=pa_l',\ l=0,\ldots,r-1,\qquad
a_l=a_l',\ l=r,\ldots,N.
\end{multline*}

The collections $\overline{d}'=\left(d_1',\ldots,d_k'\right)$ and $\overline{a}'=\left(a_0',\ldots,a_N'\right)$ form a Cartier regular pair.
Indeed, since $(\overline{d},\overline{a})$ is a Cartier regular pair, for any $l$ and $u$ the weight $a_l$ divides
$d_u$. Thus, if $l<r$, then $a_l'=a_l/p$ divides $d_u'=d_u/p$,
and if $l\ge r$, then $a_l'=a_l$ also divides $d_u'=d_u/p$, since $p$-adic valuation of
$d_u$ is greater then or equal to $s$.
This shows that $(\overline{d}',\overline{a}')$ is a Cartier pair.
It is obviously regular as well, since any $k+1$ different numbers $a_l'$
are coprime.
If $(\overline{d}',\overline{a}')$ is a Calabi--Yau pair, put
$$
\beta_l=1,\ l=0,\ldots,r-1,\qquad
\beta_l=p,\ l=r,\ldots,N.
$$
One can see that this gives the assertion of the proposition in this case.
If $(\overline{d}',\overline{a}')$ is of general type, then
%
by the induction assertion there exist positive numbers $\beta_0',\ldots,\beta_N'$
such that
$$
d_1'+\ldots+d_k'=\beta_0'a_0'+\ldots+\beta_N'a_N'.
$$
Put
\begin{equation*}
\label{equation: dividing by p}
\beta_l=\beta_l',\ l=0,\ldots,r-1,\qquad
\beta_l=p\beta_l',\ l=r,\ldots,N.
\end{equation*}
Then~\eqref{equation:splitting} holds, which proves the proposition.

Now assume that $(\overline{d}',\overline{a}')$ is Fano.
If $s>1$, then $a_l'>1$ for all $l$, so by Corollary~\ref{corollary:no 1} the pair $(\overline{d}',\overline{a}')$ is of general type.
Thus, we can assume that $s=1$ for all primes that divide $a_l$ for some $l$.
In other words, no square of a prime number divides any $a_l$.
Let $p,p_1,\ldots, p_v$ be prime divisors of all $a_l$.
We may assume that $v\ge 1$, that is there are at least two different
prime divisors of the numbers $a_l$, otherwise, since $a_l>1$ for all $l$,
$p$ divides all $a_l$, and
the pair $(\overline{d},\overline{l})$ is not regular.
In particular, we can choose $p$ such that $p>2$.
Thus for all $u=1,\ldots,k$ we have $d_u=\alpha_u\cdot p\cdot p_1\cdot\ldots\cdot p_v$ for some positive integers $\alpha_u$.
We can assume that $\alpha_u=1$ for all $u$.
Indeed, if we replace $d_u$ by $d_u/\alpha_u$,
we get a Cartier regular pair which is of general type by Corollary~\ref{corollary:no 1}.
If the assertion of the proposition holds for it, then
it holds for the initial regular pair, since if
$$
d_1/\alpha_1+\ldots+d_k/\alpha_k=\beta_0 a_0+\ldots+\beta_N a_N,
$$
then
$$
d_1+\ldots+d_k=\left(\beta_0+\sum_{u=1}^k\left(\frac{(\alpha_u-1)d_u}{\alpha_u a_0}\right)\right) a_0 +\beta_1 a_1+\ldots+\beta_N a_N.
$$

Denote $e=|\{l\mid a_l'=1\}|$ and $f=|\{l\mid  a_l'=d_u'\mbox{ for some $u$}\}|$.

Let 
$f=k$, so that 
there exist $k$ different numbers $r_l$ such that $a_{r_l}'=d'_l=p_1\cdot\ldots\cdot p_v$.
By Remark~\ref{remark: regular linear cone} applied to $(\overline{d}',\overline{a}')$, we get
$$
\overline{d}=((pc)^k),\qquad \overline{a}=(c^k,p^r)
$$
for $c=p_1\cdot\ldots\cdot p_v>1$ coprime with $p$, where 
$k\ge r=N+1-k>1$. 
If $k=2$, then $r=2$, so that $\overline{d}=(pc,pc)$ and $\overline{a}=(c,c,p,p)$. The assertion of
the proposition in this case holds by
\begin{equation}
\label{equation:two degrees} \tag{$\bigoplus$}
pc+pc=1\cdot c+(p-1)\cdot c+1\cdot p+(c-1)\cdot p.
\end{equation}
Let $k>2$.
If $r=2$, then one can set
$$
\beta_0=1,\ \beta_1=p-1,\ \beta_2=p,\ \ldots,\beta _{k-1}=p,\
\beta_{k}=1,\ \beta_{k+1}=c-1.
$$
If $r=3$ and $c>2$, then one can set
$$
\beta_0=1,\ \beta_1=p-1,\ \beta_2=p,\ \ldots,\beta _{k-1}=p,\
\beta_{k}=1,\ \beta_{k+1}=1,\ \beta_{k+2}=c-2.
$$
If $r=3$ and $c=2$, then one can set
$$
\beta_0=1+(k-2)\cdot p-k,\ \beta_1=1,\ \ldots,\beta _{k-1}=1,\
\beta_{k}=2,\ \beta_{k+1}=1,\ \beta_{k+2}=1.
$$
Let $r\ge 4$.
Using~\eqref{equation:two degrees}, one can delete degrees $pc$, $pc$ and weights $c$, $c$, $p$, $p$
from $(\overline{d},\overline{a})$ and get a Cartier regular pair. By Corollary~\ref{corollary:no 1}, it is of general type.
Moreover,
$$
(N-4)-(k-2)=N-k-2=r-3>0.
$$ 
Thus, the assertion of the proposition in this case holds by induction.

Now let $f<k$.
By Theorem~\ref{corollary:PST}, one has $e>k-f>0$, so that $e\ge 2$.
Moreover, $f>0$, otherwise by Theorem~\ref{corollary:PST} one has $e>k$, while $e$ is the number of weights $a_l$ equal to $p$, so that
$e\le k$.
Let $f=1$. Then $k\ge e>k-1$, so that $e=k$. By Theorem~\ref{corollary:PST},
this can happen only when
$$
\overline{d}=\left((6p)^k\right),\qquad \overline{a}=\left(6,2^k,3^k,p^k\right);
$$
however this pair is not regular.

Now assume that $f>1$.
We have $d_1=\ldots=d_k=p\cdot p_1\cdot\ldots\cdot p_v$,
two weights, say $a_0$ and $a_1$, are equal to $p$, and
two weights, say $a_{2}$ and $a_3$, are equal to $p_1\cdot\ldots\cdot p_v$.
Let us delete two degrees and four weights $a_0$, $a_1$, $a_{2}$, and $a_3$ from the regular pair and get the pair $(\overline{d}_0,\overline{a}_0)$.
It is a Cartier regular pair of general type (again by Corollary~\ref{corollary:no 1}).

Assume that $N=k+1$.
Then
$\overline{a}=(p,p,p_1\cdot\ldots\cdot p_v,p_1\cdot\ldots\cdot p_v,a_4,\ldots,a_{k+1})$,
so we can take
$$
\beta_0=1,\ \beta_1=(p_1\cdot\ldots\cdot p_v-1),\ \beta_2=1,\ \beta_3=p-1, \beta_l=\frac{p\cdot p_1\cdot\ldots\cdot p_v}{a_l},
\ l=4,\ldots,k+1,
$$
to prove the proposition.
Assume that $N=k+2$.
Since $a_0=a_1=p$, we have $k>1$. If $k=2$, then
$$
\overline{d}=(p\cdot p_1\cdot\ldots\cdot p_v,p\cdot p_1\cdot\ldots\cdot p_v),
\ \overline{a}=(p,p,p_1\cdot\ldots\cdot p_v,p_1\cdot\ldots\cdot p_v,a_4),
$$
so the pair $(\overline{d},\overline{a})$ is not regular, because $a_4$ divides $p\cdot p_1\cdot\ldots\cdot p_v$.
If $k>2$, then
$$
\overline{a}=(p,p,p_1\cdot\ldots\cdot p_v,p_1\cdot\ldots\cdot p_v,a_4,a_5\ldots,a_{k+2}),
$$
so we can take
\begin{multline*}
\beta_0=1,\ \beta_1=(p_1\cdot\ldots\cdot p_v-1),\ \beta_2=1,\ \beta_3=p-1,
\ \beta_4=\frac{p_1\cdot\ldots\cdot p_v}{a_4}, \\
\ \beta_5=\frac{(p-1)\cdot p_1\cdot\ldots\cdot p_v}{a_5},
\ \beta_l=\frac{p\cdot p_1\cdot\ldots\cdot p_v}{a_l},
\ l=6,\ldots,k+2
\end{multline*}
to prove the proposition.

Finally, assume that $N-k>2$, so that $(N-4)-(k-2)=N-k-2>0$.
By induction, one can find the splitting analogous to~\eqref{equation:splitting} for $(\overline{d}_0,\overline{a}_0)$.
Together with
$$
p\cdot p_1\cdot\ldots\cdot p_v+p\cdot p_1\cdot\ldots\cdot p_v=
1\cdot p+(p_1\cdot\ldots\cdot p_v-1)\cdot p+1\cdot \left(p_1\cdot\ldots\cdot p_v\right)+(p-1)\cdot \left(p_1\cdot\ldots\cdot p_v\right)
$$
it gives the required splitting for $(\overline{d},\overline{a})$.
\end{proof}

\begin{remark}
\label{remark:points}
Note that Proposition~\ref{proposition:h0n Cartier general type} does not hold for $N=k$. Indeed, let $p_0,\ldots,p_{N}$ be
different primes, let $P=p_0\cdot\ldots\cdot p_N$, and let $P_l=\frac{P}{p_l}$, so that $P_l$ is a product of all
primes except for $p_l$. Consider the regular pair $\left((P^N),(P_0,\ldots, P_{N})\right)$. Obviously, it
is Cartier and of general type. Assume that
$$
NP=\beta_0\cdot P_0+\ldots+\beta_N \cdot P_N,
$$
such that $\beta_l>0$.
Then, since $p_m$ divides $P,P_0,\ldots,P_{m-1},P_{m+1},\ldots, P_N$, it also divides $\beta_m$.
Thus,
$$
\beta_0\cdot P_0+\ldots+\beta_N \cdot P_N\ge p_0\cdot P_0+\ldots+p_N\cdot P_N=(N+1)\cdot P>NP,
$$
which is a contradiction. Note that the complete intersection of $N$ hypersurfaces
of degrees $P$ in $\PP(P_0,\ldots,P_N)$ is not well formed; after wellformization
it becomes a complete intersection of $N$ hyperplanes in $\PP^N$, that is just a point.
\end{remark}

Now we are ready to prove Theorem~\ref{theorem:maximal Hodge level}.

\begin{proof}[Proof of Theorem~\ref{theorem:maximal Hodge level}]
Assertions (i) and (ii) obviously follow from Theorem~\ref{theorem:middle-Hodge-numbers} and Corollary~\ref{corollary:h0n}.
Prove the assertion (iii).

Let $X$ be a smooth well formed weighted complete intersection of Cartier divisors, and let it be of general type. By Proposition~\ref{proposition:h0n Cartier general type}, together with Propositions~\ref{proposition:smooth-condition} and~\ref{proposition:Cartier},
there exist positive integers $\beta_0,\ldots,\beta_N$ such that
$$
d_1+\ldots+d_k=\beta_0\cdot a_0+\ldots+\beta_N\cdot a_N.
$$
Now 
the assertion (iii) follows from Corollary~\ref{corollary:sum decomposition}.
\end{proof}

\begin{problem}
\label{problem:general type maxmal}
Generalize the assertion (iii) of Theorem~\ref{theorem:maximal Hodge level} to smooth well formed weighted complete intersections
of general type (not only intersections of Cartier divisors).
\end{problem}

This can be done for codimension at most $2$ by the following proposition, shared with us by C.\,Shramov.
To prove it we need the two-coins case of Frobenius Coin Problem.

\begin{theorem}[\cite{Sy82}]
\label{theorem:coins}
Let $a_0$, $a_1$ be two coprime positive integers. Then for any integer $\mbox{$m>a_0a_1-a_0-a_1$}$ there exist non-negative integers $\beta_0$, $\beta_1$
such that $m=\beta_0\cdot a_0+\beta_1\cdot a_1$.
\end{theorem}

Now solve Problem~\ref{problem:general type maxmal} for the case of codimension $1$ and $2$.

\begin{proposition}
\label{proposition:codim 2}
The following holds.
\begin{itemize}
\item[(i)]
Let $(\overline{d},\overline{a})$, where $\overline{d}=(d_1)$ and $\overline{a}=(a_0,\ldots,a_N)$, $N>1$, be a regular pair of general type.
Then there exist positive integers $\beta_0,\ldots, \beta_N$ such that
\begin{equation*}
\label{equation:splitting codim 1}
d_1=\beta_0\cdot a_0+\ldots+\beta_N\cdot a_N.
\end{equation*}
\item[(ii)]
Let $(\overline{d},\overline{a})$, where $\overline{d}=(d_1,d_2)$ and $\overline{a}=(a_0,\ldots,a_N)$, $N>2$, be a regular pair of general type.
Assume that there exist non-negative integers $\gamma_0,\ldots,\gamma_N$ and $\mu_0,\ldots,\mu_N$ such that
$d_1=\sum \gamma_l a_l$ and $d_2=\sum \mu_l a_l$.
Then there exist positive integers $\beta_0,\ldots, \beta_N$ such that
\begin{equation*}
\label{equation:splitting codim 2}
d_1+d_2=\beta_0\cdot a_0+\ldots+\beta_N\cdot a_N.
\end{equation*}
\end{itemize}
\end{proposition}

\begin{proof}
%
%

In the assertion (i) setup the regular pair $(\overline{d},\overline{a})$ is Cartier, so it follows from Proposition~\ref{proposition:h0n Cartier general type}.
Let us prove the assertion (ii).
We follow the proof of Proposition~\ref{proposition:h0n Cartier general type}.
That is, let us prove the proposition by induction on $d_1+d_2$.
The base (when $\overline{d}=(2,3)$ and $\overline{a}=(1,1,1,1)$) is trivial.
Denote $i=d_1+d_2-\sum a_l>0$.
We can assume that $a_l>1$ for all $l$, otherwise we can set
 $\beta_m=1$, $m\neq l$, and $\beta_l=i+1$.
Moreover, we can assume that $a_l\neq d_u$ for all $l$ and $u$.
Indeed, if, say, $d_2=a_N$, then the pair
$((d_1),(a_0,\ldots,a_{N-1}))$ is regular and of general type, so by assertion (i) there exist positive numbers
$\beta_0,\ldots,\beta_{N-1}$ such that
$$
d_1=\beta_0\cdot a_0+\ldots+\beta_{N-1}\cdot a_{N-1}.
$$
Then
$$
d_1+d_2=\beta_0\cdot a_0+\ldots+\beta_{N-1}\cdot a_{N-1}+1\cdot a_N.
$$

Thus, assume that $a_l\neq 1$ and $a_l\neq d_1,d_2$ for $l=0\ldots,N$.
Let all $a_l$ be pairwise coprime.
Assume that all $a_l$ divide $d_1$.
Then, by assertion (i), there exist positive integers $\beta_0,\ldots,\beta_N$ such that
$d_1=\sum \beta_l\cdot a_l$. This means that
$$
d_1+d_2=(\beta_0+\mu_0)\cdot a_0+\ldots+(\beta_N+\mu_N)\cdot a_N,
$$
so the assertion (ii) follows. 
Assume that $a_l$ does not divide $d_1$ for some $l$.
Up to permutation we can assume that $a_0,\ldots, a_{r-1}$, $r\le N$, divide $d_1$, while $a_r,\ldots,a_N$ are not.
Note that in this case the pairs $((d_1),(a_0,\ldots,a_{r-1}))$ and $((d_2),(a_r,\ldots,a_{N}))$ are regular.
If $r=1$, then $d_1=\alpha a_0$, and by assertion (i) there exist positive integers $\beta_1,\ldots, \beta_N$
such that
$$
d_2=\beta_1\cdot a_1+\ldots+\beta_N\cdot a_N,
$$
so that
$$
d_1+d_2=\alpha\cdot a_0+\beta_1\cdot a_1+\ldots+\beta_N\cdot a_N.
$$
If $r=2$, then, since $d_2$ is divisible by coprime numbers $a_r,\ldots,a_N$ and $\mbox{$N-r+1=N-1>1$}$,
one has $d_2-a_r-\ldots-a_N>0$.
Moreover, $d_1-a_0-a_1\ge a_0a_1-a_0-a_1$.
By Theorem~\ref{theorem:coins}, there exist non-negative integers $\beta_0$, $\beta_1$ such that
$$
(d_1-a_0-a_1)+(d_2-a_2-\ldots-a_N)=\beta_0\cdot a_0+\beta_1\cdot a_1,
$$
so that
$$
d_1+d_2=(\beta_0+1)\cdot a_0+(\beta_1+1)\cdot a_1+1\cdot a_2+\ldots+1\cdot a_N.
$$
Thus, we may assume that $r>2$ and, similarly, $N+1-r>2$.
By assertion (i) there exist positive integers $\beta_0,\ldots,\beta_N$ such that
$$
d_1=\beta_0\cdot a_0+\ldots+\beta_{r-1}\cdot a_{r-1},\qquad d_2=\beta_r\cdot a_r+\ldots+\beta_{N}\cdot a_{N}.
$$
This means that
$$
d_1+d_2=\beta_0\cdot a_0+\ldots+\beta_{r-1}\cdot a_{r-1}+\beta_r\cdot a_r+\ldots+\beta_{N}\cdot a_{N}.
$$

Finally, assume that there is a prime $p$ such that there exist two weights $a_l$ and $a_m$ divisible by $p$.
Then $d_1$ and $d_2$ are divisible by $p$ as well. Consider the pair $(\overline{d}',\overline{a}')$,
where $d_1'=\frac{d_1}{p}$, $d_2'=\frac{d_2}{p}$, $a_l'=\frac{a_l}{p}$, $a_m'=\frac{a_m}{p}$, and $a_w'=a_w$ for all $w\neq l,m$.
Note that it is regular, because $p$ does not divide $a_s$ for $s\neq l,m$.
If $(\overline{d}',\overline{a}')$ is Calabi--Yau or of general type, then similarly to the proof of Proposition~\ref{proposition:h0n Cartier general type},
one can derive assertion (ii) directly or by induction.
Assume that $(\overline{d}',\overline{a}')$ is Fano.
By Theorem~\ref{corollary:PST}, this is possible
only if either $d_1'=a_s'$ and $d_2'=a_t'$ for some different indices $s$ and $t$,
or $d_1'=a_s'$ for some $s$ and $d_2'\neq a_t'$ for all $t\neq s$,
or, similarly, $d_2'=a_s'$ for some $s$ and $d_1'\neq a_t'$ for all $t\neq s$.

Let $d_1'=a_s'$ and $d_2'=a_t'$ for some different indices $s$ and $t$.
By Remark~\ref{remark: regular linear cone} applied to $(\overline{d}',\overline{a}')$ and Theorem~\ref{corollary:PST},
$(\overline{d},\overline{a})=((p\alpha,p\gamma),\left(\alpha,\gamma,p,p\right)).$
Since $p$ and $\gamma$ are coprime, we have $p\gamma-\gamma-p>0$.
Since $p$ and $\alpha$ are copime, Theorem~\ref{theorem:coins} implies that there exist non-negative integers $\beta_0$ and $\beta_2$ such that
$$
(p\alpha-\alpha-p)+(p\gamma-\gamma-p)=\beta_0\cdot \alpha+\beta_2\cdot p,
$$
so that
$$
d_1+d_2=p\alpha+p\gamma=(\beta_0+1)\cdot \alpha+1\cdot \gamma+(\beta_2+1)\cdot p+1\cdot p.
$$


Thus, without loss of generality we may assume that $d_1'=a_s'$ for some $s$ and $d_2'\neq a_t'$ for all $t\neq s$.
In this case by Theorem~\ref{corollary:PST} one has $a_l=a_m=1$ and
$(\overline{d}',\overline{a}')=((\alpha,6),(\alpha,2,3,1,1))$.
This means that $(\overline{d},\overline{a})=((\alpha p,6p),(\alpha,2,3,p,p))$.
The assertion (ii) in this case is given by
\begin{equation*}
\alpha p+6p=p\cdot \alpha+\frac{p-3}{2}\cdot 2 +1\cdot 3 +1\cdot p+ 4\cdot p. \qedhere
\end{equation*}
\end{proof}

\begin{corollary}
\label{corollary:codim 2 wci}
Let $X\subset \PP(a_0,\ldots,a_N)$ be 
      a smooth well formed weighted complete intersection of two hypersurfaces, let $\dim (X)=N-2>0$, and let it be of general type, that is 
      $i_X>0$. Then $\hh(X)=n$.
\end{corollary}

\begin{proof}
Similar to the proof of Theorem~\ref{theorem:maximal Hodge level}.
\end{proof}

\section{Hodge level for quasi-smooth weighted complete intersections of general type}
\label{section:quasi-smooth}

In this section we give an example of quasi-smooth weighted complete intersection of general type
whose Hodge level is not maximally possible.
(In fact this weighted complete intersection is either hypersurface or has codimension $2$,
depending on the parity of the dimension.)
For this we need number-theoretical results (Lemma~\ref{lemma:2n_primes}
and Proposition~\ref{proposition:prime sequence})
whose proofs are suggested to us by
M.\,Korolev.


\begin{lemma}
\label{lemma:2n_primes}
Let $n$ be a positive integer number and let $x\ge 2^n$ be a real number. Then
\begin{itemize}
  \item[(i)] if $n\ge 5$, then there are at least $n+1$ primes $p$ such that $x<p< 2x$;
  \item[(ii)] if $n\ge 7$, then there are at least $n+1$ primes $p$ such that $\frac{2}{3}x<p< x$.
\end{itemize}
\end{lemma}

\begin{proof}
Assertion (i) obviously follows from the assertion (ii).
Indeed, it can directly checked for $n=5$, so we can assume that $n\ge 6$.
Applying assertion (ii) for $2x\ge 2^{n+1}$, one sees that there are at least $n+2$ primes $p$ such that
$\frac{2}{3}\cdot (2x)<p< 2x$. In particular, there are at least $n+1$ primes
such that $x<p<2x$.

Let us prove assertion (ii).
One can check by hand that it holds for $8\le n \le 10$.
Thus we may assume that $n\ge 11$.
Let $\pi(z)$ be the number of primes less or equal then $z\in \RR$.
By~\cite[Corollary 1 of Theorem 2]{RS62} one has
$$
\frac{\alpha x}{\ln x}>\pi(x)>\frac{x}{\ln x}
$$
for $\alpha=1.25506$.
We have
\begin{multline*}
\frac{\pi(x)-\pi\left(\frac{2}{3}x\right)}{n+1}>\frac{\frac{x}{\ln x}-\frac{2\alpha}{3} \frac{x}{\ln x-\ln \frac{3}{2}}}{n+1}=
\frac{x}{\ln x\cdot (n+1)}\cdot \frac{\left(1-\frac{2\alpha}{3}\right)\cdot \ln x-\ln\frac{3}{2}}{\ln x-\ln\frac{3}{2}}\ge\\
\ge\frac{2^n}{\ln 2\cdot n \cdot (n+1)}\cdot \frac{\left(1-\frac{2\alpha}{3}\right)\cdot n\cdot \ln 2-\ln\frac{3}{2}}{\ln 2\cdot n-\ln\frac{3}{2}},
\end{multline*}
since the functions $\frac{x}{\ln x}$ and $\frac{\left(1-\frac{2\alpha}{3}\right)\cdot \ln x-\ln\frac{3}{2}}{\ln x-\ln\frac{3}{2}}$  monotonically increase if $x\ge e$.
The function
$$
H(t)=\frac{2^t}{\ln 2 \cdot t\cdot (t+1)}\cdot \frac{\left(1-\frac{2\alpha}{3}\right) \cdot t \cdot \ln 2-\ln \frac{3}{2}}{t\cdot \ln 2-\ln\frac{3}{2}}
$$
increases if $t>\frac{2}{\ln 2}$.
Since $H(10)>1$, we have
$$
\frac{\pi(x)-\pi\left(\frac{2}{3}x\right)}{n+1}>1,
$$
which gives the assertion (ii) of the lemma.
\end{proof}

Let $\delta(n)$ be the minimal non-negative number of the form
$$
1-\frac{1}{p_1}-\ldots-\frac{1}{p_n},
$$
where $p_i$ are different primes.
This number is well defined, because
one can see that $\delta(n)>0$ and that $\delta(n)$ is attained by some set of primes.
Indeed, if
$$
\frac{1}{p_1}+\ldots+\frac{1}{p_n}=1,
$$
then for $P_i=p_1\cdot\ldots\cdot p_n/p_i$ one has
$$
P_1+\ldots+P_n=p_1\cdot\ldots \cdot p_n,
$$
so that for any $i$ all numbers but $P_i$ in the latter equality is divisible by $p_i$, which is impossible.
The fact that $\delta(n)$ is attained on the finite set of collections of primes can be proved by induction.
Indeed, for $n=1$ this is trivial.
Let it hold for $n$ but does not hold for $n+1$.
This means that there exist a non-negative number $0\le a<\delta(n)$ such that
for any $\varepsilon>0$ there exist an infinite number of collections of different primes $p_1<\ldots<p_{n+1}$
such that
$$
a<1-\frac{1}{p_1}-\ldots-\frac{1}{p_{n+1}}<a+\varepsilon.
$$
This means that there exist such collection for $\varepsilon=\frac{\delta(n)-a}{2}$; moreover, we may assume that
$p_{n+1}>\frac{2}{\delta(n)-a}$ (since for any number $A$ there is only finite number of collections of primes that are
smaller then $A$).
Thus
$$
1-\frac{1}{p_1}-\ldots-\frac{1}{p_{n}}<a+\frac{\delta(n)-a}{2}+\frac{\delta(n)-a}{2}=\delta(n),
$$
which gives a contradiction.

%
%

\begin{proposition}
\label{proposition:prime sequence}
For any $n\ge 2$ there exist primes $p_1<\ldots <p_n<p_{n+1}$ such that
$$
\frac{1}{p_1}+\ldots\frac{1}{p_n}<1<\frac{1}{p_1}+\ldots\frac{1}{p_n}+\frac{1}{p_{n+1}}.
$$
\end{proposition}


\begin{proof}
The cases when $n<8$ can be checked by hand; thus we may assume that $n\ge 8$.

We claim that $\delta(n)\le \frac{1}{2^n}$.
Indeed, let us prove this by induction. One can easily check that this holds for $n=5$.
Let it hold for $m\ge 5$.
Let $p_1,\ldots,p_m$ be primes on which $\delta(m)$ is attained.
Then by Lemma~\ref{lemma:2n_primes}(i)
there are at least $m+1$ primes between $\frac{1}{\delta(m)}$ and $\frac{2}{\delta(m)}$.
Let $p_{m+1}$ be one of them different from $p_1,\ldots,p_m$.
Then
$$
\delta(m+1)\le 1-\frac{1}{p_1}-\ldots-\frac{1}{p_m}-\frac{1}{p_{m+1}}=\delta(m)-\frac{1}{p_{m+1}}
\le \delta(m)-\frac{\delta(m)}{2}=\frac{\delta(m)}{2}\le \frac{1}{2^{m+1}}.
$$

Now let 
$p_1<\ldots<p_n$ be primes such that
$$
\delta(n)=1-\frac{1}{p_1}-\ldots-\frac{1}{p_n}.
$$

Assume that $\delta(n)\ge\frac{1}{2p_n}$. Then, since 
$\frac{1}{2^n}\ge \delta(n)$, we have $p_n>2^{n-1}$.
By Lemma~\ref{lemma:2n_primes}(ii), there are at least
$n$ primes between $\frac{2}{3}p_n$ and $p_n$.
Let $p$ be one of them different from $p_1,\ldots,p_{n-1}$.
Define
$$\delta=1-\frac{1}{p_1}-\ldots-\frac{1}{p_{n-1}}-\frac{1}{p}.
$$
Since $\delta=\delta(n)+\frac{1}{p_n}-\frac{1}{p}$,
then
$$
\delta \ge \frac{1}{2p_n}+\frac{1}{p_n}-\frac{1}{p}=\frac{3/2p-p_n}{pp_n}>0,
$$
because $p>\frac{2}{3}p_n$.
However $\delta<\delta(n)$, since $p<p_n$, which contradicts the minimality of $\delta(n)$.
This shows that $\delta(n)<\frac{1}{2p_n}$.

Summarizing, we got two inequalities: $p_n<\frac{1}{2\delta(n)}$ and $\frac{1}{\delta(n)}\ge 2^n$.
The last one by Lemma~\ref{lemma:2n_primes}(i) implies that there exists a prime $p_{n+1}$ such that $\frac{1}{2\delta(n)}<p_{n+1}<\frac{1}{\delta(n)}$.
We have
$$
p_1<\ldots<p_n<\frac{1}{2\delta(n)}<p_{n+1}.
$$
The assertion of the proposition follows from
\begin{equation*}
\frac{1}{p_{n+1}}>\delta(n)=1-\frac{1}{p_1}-\ldots-\frac{1}{p_n}>0. \qedhere
\end{equation*}
\end{proof}

Now we are ready to prove the main result of this section.

\begin{proposition}
\label{proposition:quasi-smooth hypersurface}
For any $n>2$ there exists an $n$-dimensional quasi-smooth well formed weighted complete intersection $X$
of Cartier divisors such that $X$ is
of general type and $h^{0,n}(X)=0$.
\end{proposition}

\begin{proof}
Define the integer number $m>1$ by $n=2m$ if $n$ is even and $n=2m-1$ if $n$ is odd.
By Proposition~\ref{proposition:prime sequence}, there exist $m+2$ primes
$p_0<\ldots< p_{m}<p_{m+1}$ such that
$$
\frac{1}{p_0}+\ldots+\frac{1}{p_{m}}<1<\frac{1}{p_0}+\ldots+\frac{1}{p_{m}}+\frac{1}{p_{m+1}}.
$$
Put
$$
a_s=p_0\cdot\ldots\cdot p_{s-1}\cdot p_{s+1}\cdot\ldots\cdot p_{m},
$$
$s=0,\ldots,m$.
If $n$ is even, let $X$
be a general hypersurface of degree $2\cdot p_0\cdot\ldots\cdot p_{m}$ in $\PP(a_0,a_0,a_1,a_1,\ldots,a_{m},a_{m})$;
otherwise let $X$
be a general intersection of two hypersurfaces of degree $p_0\cdot\ldots\cdot p_{m}$ in $\PP(a_0,a_0,a_1,a_1,\ldots,a_{m},a_{m})$.
(On the contrary of Remark~\ref{remark:points}, it is well formed.)
Note that in both cases $X$ is quasi-smooth and well formed, and $\dim (X)=n$.
Moreover,
$$
i_X=2\cdot p_0\cdot\ldots\cdot p_{m}-2\cdot a_0-\ldots-2\cdot a_{m}=
2\cdot p_0\cdot\ldots\cdot p_{m}\cdot\left(1-\frac{1}{p_0}-\ldots-\frac{1}{p_{m}}\right)>0,
$$
so that $X$ is of general type.
On the other hand, for any $s,t$ one has
\begin{multline*}
i_X-a_s-a_t=p_0\cdot\ldots\cdot p_{m}\cdot\left(2-2\frac{1}{p_0}-\ldots-2\frac{1}{p_{m}}-\frac{1}{p_s}-\frac{1}{p_t}\right)< \\
< 2\cdot p_0\cdot\ldots\cdot p_{m}\cdot\left(1-\frac{1}{p_0}-\ldots-\frac{1}{p_{m}}-\frac{1}{p_{m+1}}\right)<0,
\end{multline*}
so there is no monomial of weight $i_X$ whose usual degree
in variables with weights $a_0,\ldots,a_{m}$ is $2$ or greater.
Finally, if for some $l$ one has $i_X=a_l$, then $i_X+2a_l=3a_l$.
However this is impossible unless $p_l=3$, since the left side of this equality is divisible by $p_l$,
while the right side is not.
For $p_l=3$ we have
\begin{multline*}
i_X-a_l=
 2\cdot p_0\cdot\ldots\cdot p_{m}\cdot\left(1-\frac{1}{p_0}-\ldots-\frac{1}{p_{m}}-\frac{1}{6}\right)<\\
< 2\cdot p_0\cdot\ldots\cdot p_{m}\cdot\left(1-\frac{1}{p_0}-\ldots-\frac{1}{p_{m}}-\frac{1}{p_{m+1}}\right)<0,
\end{multline*}
since $p_{m+1}>6$.

Thus, there are no linear monomials of weighted degree $i_X$.
Now the assertion of proposition is given by Corollary~\ref{corollary:h0n}.
\end{proof}


\end{document}